\documentclass[pdflatex,sn-mathphys-num]{sn-jnl}

\usepackage{graphicx}%
\usepackage{multirow}%
\usepackage{amsmath,amssymb,amsfonts}%
\usepackage{amsthm}%
\usepackage{mathrsfs}%
\usepackage[title]{appendix}%
\usepackage{xcolor}%
\usepackage{textcomp}%
\usepackage{manyfoot}%
\usepackage{booktabs}%
\usepackage{algorithm}%
\usepackage{algorithmicx}%
\usepackage{algpseudocode}%
\usepackage{listings}%

\theoremstyle{thmstyleone}%
\newtheorem{theorem}{Theorem}[section]
\theoremstyle{thmstyletwo}%
\newtheorem{remark}{Remark}[section]%
\newtheorem{lemma}[theorem]{Lemma}
\newtheorem{corollary}{Corollary}[section]

\theoremstyle{thmstylethree}%
\numberwithin{equation}{section}

\begin{document}

\title[A note on the growth factor in Gaussian elimination]{A note on the growth factor in Gaussian elimination for Higham matrices}

\author[1]{\fnm{Qian-Ping} \sur{Guo}}\email{guoqianpinglei@163.com}

\author*[2]{\fnm{Xian-Ming} \sur{Gu}}\email{guxianming@live.cn}

\author*[3]{\fnm{Hou-Biao} \sur{Li}}\email{lihoubiao0189@163.com}

\affil[1]{\orgdiv{School of Mathematics and Information Science}, \orgname{Henan
University of Economics and Law}, \orgaddress{\city{Zhengzhou}, \postcode{450046}, \country{P.R. China}}}

\affil[2]{\orgdiv{School of Mathematics}, \orgname{Southwestern University of Finance and Economics},
\orgaddress{\city{Chengdu}, \postcode{611130}, \country{P.R. China}}}

\affil[3]{\orgdiv{School of Mathematics sciences}, \orgname{University of Electronic Science and Technology of
China}, \orgaddress{\city{Chengdu}, \postcode{611731}, \country{P.R. China}}}


\abstract{The Higham matrix is a complex symmetric matrix $A=B+\mathrm{i}C$,
where both matrices $B$ and $C$ are real, symmetric and positive definite and
$\mathrm{i}=\sqrt{-1}$ is the imaginary unit. For any Higham matrix $A$,
Ikramov et al. showed that the growth factor in Gaussian
elimination is less than $3$. In this paper, based on the previous results, a new bound of the growth
factor is obtained by using the maximum of the condition numbers of the
matrices $B$ and $C$ for the generalized Higham matrix $A$, which strengthens this bound to $2$ and proves Higham's conjecture.}




\maketitle

\section{Introduction}\label{sec1}
In this paper, we mainly consider the complex symmetric linear system
\begin{equation}\label{eq:1.1}
Ax=b,
\end{equation}
with the coefficient matrix $A\in\mathbb{C}^{n\times n}$ and $b \in \mathbb{C}^{n}$. Such
linear systems (\ref{eq:1.1}) arise from different physical
applications, for example, modeling electromagnetic waves under
the assumption of time-harmonic variation in electromagnetic
fields, Maxwell (or Helmholtz) equations utilized a complex shift \cite{Erlangga04,Gu2015}. Moreover, the complex-valued linear system
(\ref{eq:1.1}) can be directly generated from the centered
difference discretization of the $R_{22}$-Pad\'{e} approximations in
the time integration of parabolic partial differential equations (PDEs)
\cite{OAK} and in direct frequency domain analysis of an
$n$-degree-of-freedom linear system \cite{AFF,Simoncini02}. In addition, there
are some other matrix computational problems \cite{Frommer08,Arbenz04,BEREUX05,sssx} that involve the solution of a linear system (\ref{eq:1.1}).
Therefore, the research on (sparse) direct and iterative solvers for the linear system (\ref{eq:1.1})
is highly needed.

Next, for convenience, let $M_{n}(\mathbb{C})$ denote the set of
$n\times n$ complex matrices and $A$ be a nonsingular matrix in
$M_{n}(\mathbb{C})$. $\lambda_{1}$ and $\lambda_{n}$ are the largest
and smallest eigenvalues of $A^{\ast} A$, respectively.

Recently, a \emph{complex symmetric positive definite} (CSPD) matrix
arising from the linear system (\ref{eq:1.1})
in Gaussian elimination without pivoting was firstly studied by
Higham in \cite{HN}, which is called by \emph{\textbf{Higham
matrices}} in \cite{A.K}. Subsequently, the paper \cite{A.K} gave a
broader class of complex matrices---\emph{\textbf{generalized Higham
matrices}} (sometimes they are also called accretive-dissipative matrices
\cite{ML}), i.e., for any $A\in M_{n}(\mathbb{C})$, if, its
Toeplitz decomposition \cite{NJ} (sometimes also called the
Hermitian decomposition \cite{ML}).
 $$
  A=B+\mathrm{i}C
 $$
satisfies that
$$
\mathrm{i}=\sqrt{-1},\ B=B^{\ast}>0,\ C=C^{\ast}>0,
 $$
where $B^{\ast}$ is the conjugate transpose of $B$, then the matrix
$A$ is said to be a generalized Higham matrix and denoted by $A \in
M_{n}^{++}(\mathbb{C})$. Here, the sign $\geq$ is usually called the
Loewner partial order of Hermitian matrices; i.e., we write $B\geq
C$ if the matrix $B-C$ is positive semidefinite, similarly, $B>C$
means that $B-C$ is positive definite. In addition, a related class
of matrices defined by
 $$
  A=B+\mathrm{i}C,  ~~~~B=B^{\ast}>0,~~~~C=C^{\ast}<0,
 $$
will be denoted by $M_{n}^{+-}(\mathbb{C})$ as in \cite{A.K}.

This paper is a continuation of \cite{A.K}, and both originate
from Higham's paper \cite{HN}. As is well known, for the complex linear systems
(\ref{eq:1.1}), the growth factor $\rho_{n}(A)$ in Gaussian
elimination for $A$ is defined by
\begin{equation}\label{eq:1.2}
\rho_{n}(A)=\frac{\max\limits_{i,j,k}|a^{(k)}_{ij}|}{\max\limits_{i,j}|a_{ij}|},
\end{equation}
where $A\triangleq(a_{ij})$, $A^{(k)}\triangleq(a_{ij}^{(k)})$ and
$A^{(k)}$ is the matrix obtained through the application of the
first $k$ steps of Gaussian elimination to the matrix $A$, see
e.g. \cite{A.K,GW}. In particular, $A^{(n-1)}$ is the upper
triangular matrix resulting from the $LU$ factorization of $A$.

Obviously, for the matrix $A$ in the linear system (\ref{eq:1.1}), if
one is able to prove a satisfactory priori bound for $\rho_{n}(A)$,
then it is safe not to pivot in computing the $LU$ factorization
of the matrix $A$ (or to choose diagonal pivots based on other
considerations such as sparsity preservation) \cite{A.K,NJ}. In addition, this insight holds significant guiding importance for the development of efficient algorithms to solve complex linear systems (\ref{eq:1.1}).

For any Higham matrix $A$, the growth factor in Gaussian elimination
\begin{equation}\label{eq:1.4}
{\rho _n}(A) \le 2
\end{equation}
was firstly conjectured by Higham \cite[p.210]{NJ}.
An incorrect proof was given in \cite{HN}, but Ikramov et al. \cite{A.K} subsequently showed that
\begin{equation}\label{eq:1.5}
 {\rho _n}\left( A \right) < 3
  \end{equation}
for any Higham matrix $A$. In addition, if the Higham matrix is
extended by allowing $B$ and $C$ to be arbitrary Hermitian positive
definite matrices, then
\begin{equation} \label{eq:1.6}
{\rho _n}\left( A \right) < 3\sqrt{2}.
\end{equation}
Moreover, Ikaramov et al. noted that the above bound (\ref{eq:1.6})
remains true when $B$ or $C$ or both are negative (rather than
positive) definite \cite{A.K}.

For a very restricted subset of Higham matrices, i.e.,
when $B=I_{n}$ and $C$ is real, symmetric and positive definite, the
authors in \cite{I.K} proved the better bound
$${\rho _n}(A) \le \frac{{1 + \sqrt {17} }}{4} \approx 1.28078
\cdots.
$$

In addition, George and Ikramov \cite{AK} assumed that $B$
and $C$ with being positive definite, satisfy the
inequality
\[
C \le \alpha B,\;\;\;\alpha  \ge 0,
\]
and they established a bound for the growth factor ${\rho
_n}(A)$ that has the limit $1$ as $\alpha  \to 0 $.

Recently, Lin \cite{MLin} proved that if $A$ is a generalized Higham matrix,
then the growth factor for such a matrix $A$ in Gaussian elimination is less than $4$. Especially when $A$ is a Higham matrix, the growth factor is less than $2\sqrt{2}$. Subsequently, Yang \cite{Yang} gave a minor improvement of the results
in Lin \cite{MLin}, which is much closer to the final solution of Higham's conjecture; refer to \cite[Theorem 4]{Yang}. This paper \cite{Drury} derives an upper bound for the growth factor via the decomposition and Fischer-type determinantal inequalities, it shows that $\rho_n(A)\leq sec^3(\alpha)$. In particular, when $\alpha=\frac{\pi}{4}$,
$$
\rho_n(A)\leq 2\sqrt{2}.
$$

However, as the authors in \cite{A.K} pointed out, in no case have they
observed a growth bigger than the Higham guess of $2$ from extensive
numerical experiments with Higham matrices. Therefore, they believed
that the bound (\ref{eq:1.4}) is correct and took the proof of this
bound for an open problem (see also Problem 10.12 in \cite{NJ}).

In this work, we continue studying this open problem and then give a
new result
$$
0\leq\frac{{4\kappa }}{{{{\left( {1 + \kappa } \right)}^2}}} \le
{\rho _n}\left( A \right) \le \frac{{2\left( {1 + {\kappa ^2}}
\right)}}{{{{\left( {1 + \kappa } \right)}^2}}} \leq 2,
$$
for the generalized Higham matrix $A$, where $\kappa\in [1,+\infty)$
is the maximum of the condition numbers of $B$ and $C$. This
directly leads to the Higham's result (\ref{eq:1.4}) for any Higham
matrix $A$, which proves the open problem. Here, for a nonsingular matrix $A$, its spectral condition
number is denoted by $\kappa (A)\triangleq \sqrt {\frac{{{\lambda
_{\max }}\left( {{A^ * }A} \right)}}{{{\lambda _{\min }}\left( {{A^
* }A} \right)}}} $, i.e., the ratio of the
largest and smallest singular value of $A$.

This paper is organized as follows. In Section \ref{sec2}, we show some new
bounds on the growth factor, based on the condition number. In
Section \ref{sec3}, the numerical example and table are given to
illustrate our results. Finally, we draw some conclusions in Section \ref{sec4}.

\section{Results}\label{sec2}

In this section, let $A\in M_{n}(\mathbb{C})$ be partitioned as
\begin{equation}\label{eq:1.3}
A \triangleq
\begin{pmatrix}
A_{11} & A_{12} \\
A_{21} & A_{22}
\end{pmatrix} =
\begin{pmatrix}
B_{11} & B_{12}  \\
B_{21} & B_{22}
\end{pmatrix} + \mathrm{i}
\begin{pmatrix}
C_{11}  & C_{12}\\
C_{21}  & C_{22}
\end{pmatrix},
\end{equation}
where $A$ is an $n\times n$ nonsingular matrix. If $A_{11}$ is
invertible, then the Schur complement of $A_{11}$ in $A$ is denoted
by $A/{A_{11}} = {A_{22}}- {A_{21}}A_{11}^{-1}A_{12}$ (see,
\cite{W.Z}).

\begin{lemma}[\cite{I.K}]\label{lem:2.1}
Let $A$ be a CSPD matrix, then $A$ is nonsingular, and any principal
submatrix of $A$ and any Schur complement in $A$ are also CSPD
matrices.
\end{lemma}
Obviously, Lemma \ref{lem:2.1} shows that, being a CSPD matrix is an
hereditary property of active submatrices in Gaussian elimination.
\begin{lemma}[\cite{I.K}]\label{lem:2.2}
The largest element of a CSPD matrix $A$ lies on its main diagonal.
\end{lemma}

Thus, for any CSPD matrix $A$, the definition (\ref{eq:1.2}) can be
replaced by
$$
{\rho _n}\left( A \right) = \frac{\max\limits_{j,k}\left|a_{jj}^{(k)}\right|}
{\max\limits_{j}|a_{jj}|},
$$
which greatly simplifies the analysis on bounding the growth factor for a CSPD matrix $A$.

\begin{lemma}[\cite{H.J}]\label{lem:3.1}
If $B$ is a nonzero $n\times n$ positive definite matrix having
eigenvalues $\lambda_{1}\geq\lambda_{2}\geq \cdots \geq
\lambda_{n}$, then for all orthogonal vectors $x,y\in
\mathbb{C}^{n}$, and $x^{\ast}$ denotes the conjugate transpose of
$x$, the following equality holds,
\begin{equation}\label{eq:3.1}
|{x^ * }By{|^2} \le {\left( {\frac{{{\lambda _1} - {\lambda
_n}}}{{{\lambda _1} + {\lambda _n}}}} \right)^2}\left( {{x^ * }Bx}
\right)\left( {{y^ * }By} \right).
\end{equation}
\end{lemma}

\begin{lemma}[\cite{FZ}]\label{Lem:3.1x}
Let $B$ be as in Lemma \ref{lem:3.1}, then for any $n\times p$ matrix $X$
satisfying ${X^ *}X = {I_p}$, where $X^{\ast}$ means the conjugate
transpose of the matrix $X$, we have that
\begin{equation}\label{eq:3.2}
{X^ * }{B^{ - 1}}X \le \frac{{{{\left( {{\lambda _1} + {\lambda _n}}
\right)}^2}}}{{4{\lambda _1}{\lambda _n}}}{\left( {{X^ * }{B}X}
\right)^{ - 1}}.
\end{equation}
\end{lemma}

By Lemma \ref{Lem:3.1x}, it is easy to obtain the
following lemma.

\begin{lemma}[\cite{FZ}]\label{lem:3.1y}
Let $B = \begin{pmatrix}
B_{11}  & B_{12} \\
B_{21}  & B_{22}
\end{pmatrix}$ be an $n\times n$ Hermitian positive definite matrix, where $B_{22}$ is any $k\times k$ principal submatrix of $B$ ($k>0$), then
\begin{equation}\label{eq:3.3}
B_{21}B_{11}^{-1}B_{12}\le \left[\frac{1 - \kappa(B)}{1
+ \kappa(B)}\right]^2B_{22},
\end{equation}
where $\kappa ({B})$ is the condition number of ${B}$.
\end{lemma}

\begin{theorem}[\cite{H.J}]\label{th:3.2}
Let $B$ be a Hermitian positive definite matrix, then
$\lambda_{n-t+i}(B)\leq\lambda_{i}(B_{t})\leq\lambda_{i}(B)$,
($i=1,2,\cdots,t$), where $B_{t}=B(i_{1},\cdots,i_{t})$ is the
$t\times t$ principal submatrix of $B$.
\end{theorem}

\begin{corollary}[\cite{Wang}]\label{cor:3.2}
Let $B$ be a Hermitian positive definite matrix, and partitioned as
in Lemma \ref{lem:3.1y}, then $\kappa(B)>\kappa(B_{11})$.
\end{corollary}

\begin{lemma}[\cite{ML}]\label{lem:2.7}
Let $A=B+\mathrm{i}C$, $B=B^{*}$, $C=C^{*}$, be partitioned as in Eq.
\eqref{eq:1.3}, if $B_{11}$, $C_{11}$ are invertible, then
\begin{equation}\label{eq:3.4}
A/A_{11}=B/B_{11}+\mathrm{i}(C/C_{11})+X\left(B_{11}^{-1}-\mathrm{i}C_{11}^{-1}\right)^{-1}X^{*},
\end{equation}
where $X=B_{21}B_{11}^{-1}-C_{21}C_{11}^{-1}$.
\end{lemma}


\begin{corollary}[\cite{ML}]\label{cor:3.3}
Let $A=B+\mathrm{i}C$, $B=B^{*}$, $C=C^{*}$ be a generalized Higham
matrix and be partitioned as in Eq. (\ref{eq:1.3}), if $A/A_{11}=R+iS$
is its Hermitian decomposition, then $R\geq B/B_{11}$, $S\geq
C/C_{11}.$
\end{corollary}
\begin{theorem}\label{th:2.8}
Let $A=B+\mathrm{i}C$, $B=B^{*}$, $C=C^{*}$ be a generalized Higham
matrix and be partitioned as in Eq. (\ref{eq:1.3}), if $A/A_{11}=R+iS$
is its Hermitian decomposition, and $B_{11}$, $C_{11}$ are invertible, then
$$
Re[X(B_{11}^{-1}-\mathrm{i}C_{11}^{-1})^{-1}X^{*}]\geq0,\quad Im[X(B_{11}^{-1}-\mathrm{i}C_{11}^{-1})^{-1}X^{*}]\geq0.
$$
\end{theorem}
\begin{proof} By applying Lemma \ref{lem:2.7} and Corollary \ref{cor:3.3}, we get
$$
R=B/B_{11}+Re[X(B_{11}^{-1}-\mathrm{i}C_{11}^{-1})^{-1}X^{*}]\geq B/B_{11},
$$
$$
S=C/C_{11}+Im[X(B_{11}^{-1}-\mathrm{i}C_{11}^{-1})^{-1}X^{*}]\geq C/C_{11},
$$
which implies $Re[X(B_{11}^{-1}-\mathrm{i}C_{11}^{-1})^{-1}X^{*}]\geq 0$ and $Im[X(B_{11}^{-1}-\mathrm{i}C_{11}^{-1})^{-1}X^{*}]\geq0$.
\end{proof}

\begin{figure}[htbp]
\includegraphics[width=2.4in,height=2.8in]{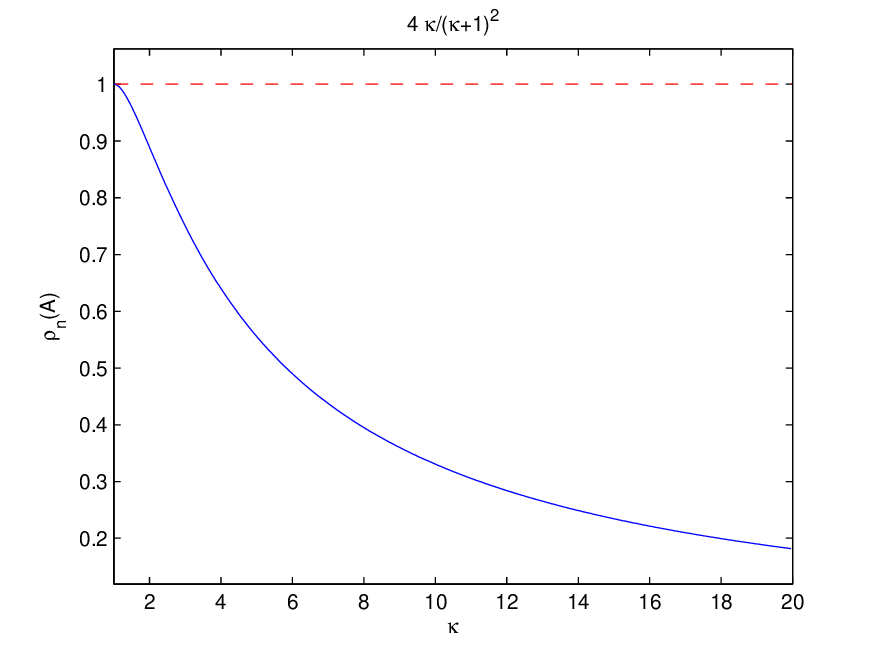}
\includegraphics[width=2.4in,height=2.8in]{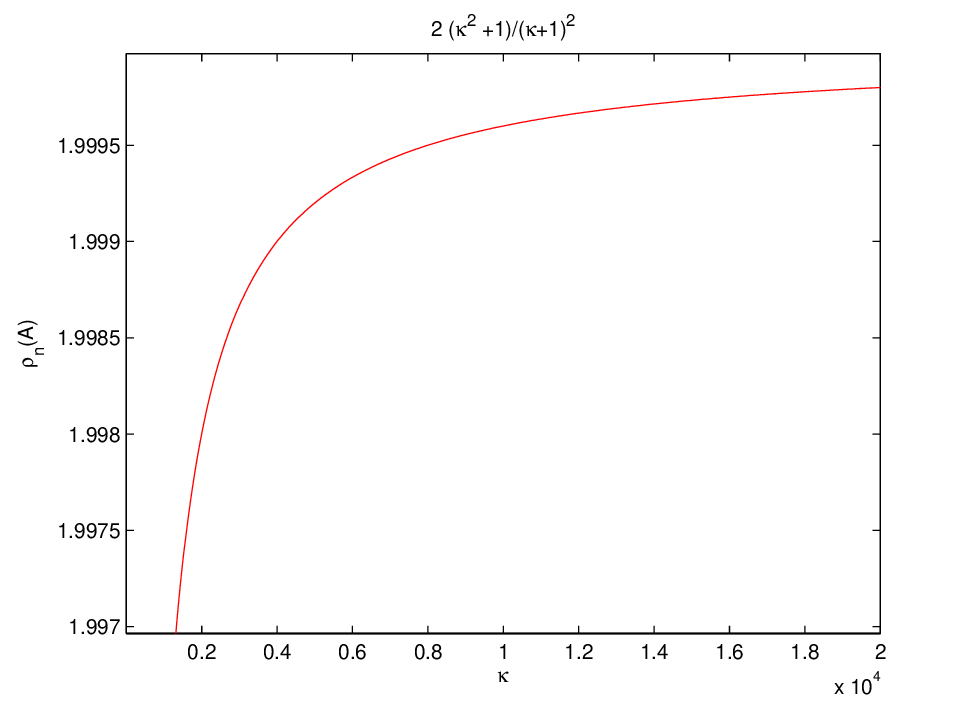}
\caption{\small Left: the variation curve of
$f(\kappa)=\frac{4\kappa}{(1+\kappa)^2}$ with $\kappa$ decreasing.
Right: the variation curve of $g(\kappa)=\frac{2(1+\kappa
^2)}{(1+\kappa )^2}$ with $\kappa$ increasing.}
\label{fig1}
\end{figure}

Next, we give our main result.
\begin{theorem}\label{th:3.3}
Let $A$ be a generalized Higham matrix, then
\begin{equation}\label{eq:3.5}
\frac{{4\kappa }}{{{{\left( {1 + \kappa } \right)}^2}}} \le {\rho
_n}\left( A \right) \le \frac{{2\left( {1 + {\kappa ^2}}
\right)}}{{{{\left( {1 + \kappa } \right)}^2}}},
\end{equation} where $\kappa$ is the maximum of the condition numbers of $B$ and $C$.
\end{theorem}
\begin{proof} Fix the integer $k\in\{1,2,\cdots,n-1\}$ and $j$,
where $j\geq k+1$. Denote $A_{k}$ by the leading principal order $k$
submatrix of $A$.

We consider the $(k+1)\times (k+1)$ matrix
$$
A_{kj} = \begin{pmatrix}
   A_{k} & \alpha   \\
   \beta^{\top} & a_{jj}  \\
\end{pmatrix} = B_{kj}+\mathrm{i}C_{kj},
$$
where
$$\alpha^{\top} = \begin{pmatrix}
   a_{1j}, a_{2j}, \cdots, a_{kj}
   \end{pmatrix}~\mathrm{and}~
\beta^{\top} = \begin{pmatrix}
a_{j1}, a_{j2}, \cdots, a_{jk}
   \end{pmatrix},
$$
$$
B_{kj} = \begin{pmatrix}
   B_{k} & b   \\
   b^{*} & b_{jj}   \\
\end{pmatrix}~\mathrm{and}~
C_{kj} = \begin{pmatrix}
   C_{k} & c   \\
   c^{*} & c_{jj}   \\
\end{pmatrix}.
$$
Note that $A_{kj}$, $B_{kj}$ and $C_{kj}$ are principal order $k+1$
submatries of $A$, $B$ and $C$, respectively.

It is easy to see that $a^{(k)}_{jj}$ can be obtained by performing
block Gaussian elimination in $A_{kj}$; i.e.,
$$
a_{jj}^{(k)} = {a_{jj}} - {\beta^\top}A_k^{ - 1}\alpha.
$$
Setting $a_{jj}^{(k)} = \eta + \mathrm{i}\gamma,~\eta,~\gamma
\in \mathbb{R}$. Since both $B_{kj}$ and $C_{kj}$ are Hermitian
positive definite, according to the result of Lemma
\ref{lem:3.1}, we have
$$
b^{*}B_k^{-1}b \le \left[\frac{1 - \kappa(B_{kj})}{1 +
\kappa(B_{kj})}\right]^2b_{jj}\ \mathrm{and}\ c^{*}C_k^{-1}c
\le \left[\frac{1 - \kappa(C_{kj})}{1 + \kappa
(C_{kj})}\right]^2c_{jj},
$$
i.e,
$$
B_{kj}/B_{k}=b_{jj}-b^*B_{k}^{-1}b\geq\frac{4\kappa(B_{kj})}{[1+\kappa(B_{kj})]^2}b_{jj},
$$
$$
C_{kj}/C_{k}=c_{jj}-c^*C_{k}^{-1}c\geq\frac{4\kappa(C_{kj})}{[1+\kappa(C_{kj})]^2}c_{jj}.
$$
Next, with the help of Corollary \ref{cor:3.3} and notice that the function
$f(x)=\frac{4x}{(1+x)^{2}}$ is decreasing on $x\in[1,+\infty)$, we
get
\begin{equation}\label{eq:3.6}
\begin{split}
|a_{jj}^{(k)}|&=|\eta+\mathrm{i}\gamma|\\
&\geq|B_{kj}/B_{k}+\mathrm{i}C_{kj}/C_{k}|\\
&\geq\sqrt{\left[\frac{4\kappa(B_{kj})}{(1+\kappa(B_{kj}))^2}b_{jj}\right]^2+\left[\frac{4\kappa(C_{kj})}{(1+\kappa(C_{kj}))^2}c_{jj}\right]^2}\\
&\geq\sqrt{\left[\frac{4\kappa_{kj}}{(1+\kappa_{kj})^2}b_{jj}\right]^2+\left[\frac{4\kappa_{kj}}{(1+\kappa_{kj})^2}c_{jj}\right]^2}\\
&=\frac{4\kappa_{kj}}{(1+\kappa_{kj})^{2}}\sqrt{b_{jj}^2+c_{jj}^2}\\
&=\frac{4\kappa_{kj}}{(1+\kappa_{kj})^{2}}|a_{jj}|.
\end{split}
\end{equation}
where $\kappa_{kj}=max\{\kappa(B_{kj}), \kappa(C_{kj})\}$.

Since $a_{jj}^{(k)}={a_{jj}}-{\beta^\top}A_k^{-1}\alpha
=b_{jj}+\mathrm{i}c_{jj}-\beta^{\top}A^{-1}_{k}\alpha$,
$B/B_{11}=b_{jj}-b^{\ast}B_{k}^{-1}b$ and $C/C_{11}=c_{jj}-c^{\ast}C_{k}^{-1}c$, by Lemma \ref{lem:2.7}, we have
\begin{equation*}
\beta^{\top}A^{-1}_{k}\alpha=b^*B_k^{-1}b + \mathrm{i}c^*C_k^{-1}c-x(B_{k}^{-1}-\mathrm{i}C_{k}^{-1})^{-1}x^{*},
\end{equation*}
where $x=b^*B_{k}^{-1}-c^*C_{k}^{-1}$. Obviously,
$$
Re(\beta^{\top}A^{-1}_{k}\alpha)=b^{\ast}B_{k}^{-1}b-Re[x(B_{k}^{-1}-iC_{k}^{-1})^{-1}x^{*}],
$$
and
$$
Im(\beta^{\top}A^{-1}_{k}\alpha)=c^{\ast}C_{k}^{-1}c-Im[x(B_{k}^{-1}-iC_{k}^{-1})^{-1}x^{*}].
$$


With the help of Theorem \ref{th:2.8}, we obtain
$$
 Re(\beta^{\top}A_{k}^{-1}\alpha)\leq b^{*}B^{-1}_{k}b\leq \left[\frac{1-\kappa(B_{kj})}{1+\kappa(B_{kj})}\right]^{2}b_{jj}
$$
and
$$
Im(\beta^{\top}A_{k}^{-1}\alpha)\leq c^{*}C^{-1}_{k}c\leq \left[\frac{1-\kappa(C_{kj})}{1+\kappa(C_{kj})}\right]^{2}c_{jj}.
$$

Since
$$
|\beta^{\top}A_{k}^{-1}\alpha|\leq\sqrt{[Re(\beta^{\top}A_{k}^{-1}\alpha)]^2+[Im(\beta^{\top}A_{k}^{-1}\alpha)]^2}
$$
and noting that $g(x)=\left(\frac{1-x}{1+x}\right)^{2}$ is increasing on $x\in[1,+\infty)$, the above inequality equivalently,
\begin{equation*}
\begin{split}
|\beta^{\top}A_{k}^{-1}\alpha|& \leq\sqrt{\left[\frac{(1-\kappa_{(B_{kj})})^2}{(1+\kappa_{(B_{kj})})^2}b_{jj}\right]^2
+\left[\frac{(1-\kappa_{(C_{kj})})^2}{(1+\kappa_{(C_{kj})})^2}c_{jj}\right]^2}\\
& \leq\left[\frac{1-\kappa_{kj}}{1+\kappa_{kj}}\right]^{2}\sqrt{b_{jj}^2+c_{jj}^2}\\
& \leq\left[\frac{1-\kappa_{kj}}{1+\kappa_{kj}}\right]^{2}|a_{jj}|.
\end{split}
\end{equation*}
Thus, we can obtain
\begin{equation}\label{eq:3.7}
\begin{split}
|a_{jj}^{(k)}| & = |a_{jj} - \beta^{\top}A_k^{-1}\alpha|\\
&\leq |{a_{jj}}| + |{\beta^\top}A_k^{ - 1}\alpha|\\
&\leq |{a_{jj}}| + \left(\frac{1-\kappa_{kj}}{1+\kappa_{kj}}\right)^{2}|a_{jj}|\\
&=\frac{2(1+\kappa_{kj}^{2})}{(1+\kappa_{kj})^{2}}|{a_{jj}}|.
\end{split}
\end{equation}
According to the above inequalities (\ref{eq:3.6}) and
(\ref{eq:3.7}), the following inequalities is obvious that
$$
\frac{4\kappa_{kj}}{(1+\kappa_{kj})^{2}}\leq\rho_{n}(A)\leq\frac{2(1+\kappa_{kj}^{2})}{(1+\kappa_{kj})^{2}}.
$$

Note that $f\left( x \right) = \frac{{4x}}{{{{\left( {1+ x}
\right)}^2}}}$ is decreasing on $x\in[1,+\infty]$, and $g\left( x
\right) = \frac{{2\left( {1 + {x^2}} \right)}}{{{{\left( {1 + x}
\right)}^2}}}$ is increasing on $x\in[1,+\infty]$ (see, Fig. \ref{fig1}), by
Corollary \ref{cor:3.2}, we have that
\begin{equation}\label{eq:3.8}
\frac{{4\kappa }}{{{{\left( {1 + \kappa } \right)}^2}}} \le {\rho
_n}\left( A \right) \le \frac{{2\left( {1 + {\kappa ^2}}
\right)}}{{{{\left( {1 + \kappa } \right)}^2}}}.
\end{equation}
The proof is complete.
\end{proof}

\begin{corollary}\label{cor:3:4}
If $A$ is an $n\times n$ generalized Higham matrix, then
\begin{equation}\label{eq:38}
0\leq\rho_{n}(A)\leq 2.
\end{equation}
\end{corollary}
\begin{proof} Since $\frac{{4\kappa }}{{{{\left({1 + \kappa }
\right)}^2}}} \le {\rho _n}\left( A \right) \le \frac{{2\left( {1 +
{\kappa ^2}} \right)}}{{{{\left( {1 + \kappa } \right)}^2}}}$, we simultaneously
calculate the limit of both sides of the inequality,
we have
$$
\mathop {\lim }\limits_{\kappa  \to \infty } \frac{{4\kappa
}}{{{{\left( {1 + \kappa } \right)}^2}}} = 0~~ \mathrm{and}~~
\mathop {\lim }\limits_{\kappa  \to \infty } \frac{{2(1 + {\kappa
^2})}}{{{{\left( {1 + \kappa } \right)}^2}}} = 2.
$$
Therefore, the result (\ref{eq:38}) is valid.
\end{proof}
\begin{remark} Obviously, the above result (\ref{eq:38}) also holds
for any Higham matrix who is just a special case of generalized Higham matrices, and hence Higham's conjecture (see
(\ref{eq:1.4})) is correct, which solves this open problem. In
addition, since $\kappa\geq 1>0$, then
$(\kappa+1)^{2}=\kappa^{2}+2\kappa+1>\kappa^{2}+1$, that is,
$\frac{2(\kappa^{2}+1)}{(\kappa+1)^{2}}<2$. So, generally speaking,
$\rho_{n}(A)< 2$, which will be shown in the following numerical experiment.
\end{remark}

\section{Numerical experiments}
\label{sec3}
In this section, a numerical example is reported to examine the effectiveness of our
proposed result. All the numerical experiments were conducted using MATLAB R2023a on 11th Gen Intel(R) Core(TM) i7-11700K @ 3.60GHz and 48.0 GB of RAM.

We consider the complex symmetric linear system
(\ref{eq:1.1}), which arises in the centered difference discretization of
the $R_{22}$-Pad\'{e} approximations in the time integration of
parabolic PDEs, further details refer to
\cite{ZZB}. For convenience, the complex coefficient symmetric matrix (see,
\cite{OAK}) may be written as
$$
A=\left(K+\frac{3-\sqrt{3}}{\tau}I\right)+\mathrm{i}\left(K+\frac{3+\sqrt{3}}{\tau}I\right),
$$
where $I$ is the identity matrix of suitable order, $\tau$ is the time step-size and
$K$ is the five-point centered difference matrix approximating the
negative Laplacian operator $L=-\Delta$ with homogeneous Dirichlet
boundary conditions, on a uniform mesh in the unit square
$[0,1]\times[0,1]$ with the mesh-size $h=\frac{1}{m+1}$.

In our tests, we take $\tau= h$. The matrix $K\in\mathbb{R}^{n\times n}$
possesses the sum of Kronecker products $K=I\otimes V_{m}+V_{m}\otimes I$,
with $V_{m}=h^{-2}\mathrm{tridiag}(-1,2,-1)\in\mathbb{R}^{m\times
m}$. Hence, $K$ is a block tridiagonal matrix, of order
$n=m^{2}$.

Denote
$$
B=K+\frac{3-\sqrt{3}}{\tau}I  ~~\mathrm{and}~~
C=K+\frac{3+\sqrt{3}}{\tau}I,
$$
and they are obviously both real symmetric positive definite. According to the above definition of the matrix $K$, it can be quickly diagonalized with its eigenvalues \cite{GDS} as follows:
\begin{equation}
\lambda_{i,j} = h^{-2}\left[4 - 2\left(\cos\frac{\pi i}{m+1} + \cos\frac{\pi j}{m+1}\right)\right],\quad i,j=1,\cdots,m.
\end{equation}
We can easily compute the eigenvalues of the matrices $B$ and $C$ as follows:
$$\lambda_{i,j}(B) = \lambda_{i,j} + \frac{3-\sqrt{3}}{\tau}\ {\rm and}\ \lambda_{i,j}(C) = \lambda_{i,j} + \frac{3+\sqrt{3}}{\tau},
$$
thus their spectral condition numbers are given as follows:
\begin{equation}
t_1 = \kappa(B) = \frac{\max\{\lambda_{i,j}(B)\}}{\min\{\lambda_{i,j}(B)\}}\ {\rm and}\
t_2 = \kappa(C) = \frac{\max\{\lambda_{i,j}(C)\}}{\min\{\lambda_{i,j}(C)\}}.
\end{equation}

\begin{table}[!htbp]
  \centering
  \caption{The changing of growth factor with the increase of condition number.}
  \label{tab1}
  \begin{tabular}{cccc}
    \toprule
    Size~($m$) & $t_1$ & $t_2$ & $L$ \\
    \midrule
    $2^5$    & 1.4186e+02 & 5.0309e+01 & 1.9722 \\
    $2^6$    & 3.3149e+02 & 1.0414e+02 & 1.9880 \\
    $2^7$    & 7.2705e+02 & 2.1219e+02 & 1.9945 \\
    $2^8$    & 1.5298e+03 & 4.2851e+02 & 1.9974 \\
    $2^9$    & 3.1423e+03 & 8.6127e+02 & 1.9987 \\
    $2^{10}$ & 6.3714e+03 & 1.7268e+03 & 1.9994 \\
    $2^{11}$ & 1.2831e+04 & 3.4580e+03 & 1.9997 \\
    \bottomrule
  \end{tabular}
\end{table}

As seen in Table \ref{tab1}, we define $L = \frac{2(1+\kappa^2)}{(1+\kappa)^2}$ and then it is evident that the results of this experiment are consistent with our
theoretical analysis.

\section{Conclusions}
\label{sec4}
In this paper, our results in Theorem \ref{th:3.3} demonstrate and answer the Higham conjecture by using the analysis of the condition number, and it was concluded that the growth factor does not exceed 2.

\section*{Acknowledgment}
The authors sincerely thank Prof.
Nicholas J. Higham for bringing his monograph \cite{NJ} to our attention and some friendly
discussions on this topic about ten years ago. However, we still fail to work out the
key proof in the past decade. In 2024, we are very sorry to hear that he passed away,
then we try to reconsider our proof, and fortunately we resolve his conjecture at last.
In addition, Dr. Minghua Lin (when he was a Ph.D. student at the University of Waterloo) provided some useful material for the manuscript, which led to a substantial improvement of
this paper.
\bibliography{sn-bibliography}

\end{document}